\documentclass[10pt,a4paper]{article} 
\usepackage{fullpage}

\newcommand{\f}{\frac}

\newcommand{\del}{\partial}

\newcommand{\im}{\operatorname{Im}}

\newcommand{\R}{\mathbb R}

\newcommand{\C}{\mathbb C}

\newcommand{\re}{\operatorname{Re}}
\newcommand{\eps}{\varepsilon}

\newcommand{\Om}{\Omega}

\newcommand{\supp}{\operatorname{supp}}
\renewcommand{\H}{\mathcal H}

\newcommand{\dist}{\operatorname{dist}}

\usepackage[english]{babel}
\usepackage{tikz}
\usetikzlibrary{arrows,decorations.pathmorphing,backgrounds,positioning,fit,petri,patterns}
\usepgflibrary{patterns}
\usepackage[utf8]{inputenc}
\usepackage{graphicx}
\usepackage[colorlinks=true,linkcolor=black,citecolor=black]{hyperref}
\usepackage{amsfonts}
\usepackage{enumerate}
\usepackage{booktabs}
\usepackage{array} 
\usepackage{paralist} 
\usepackage{subfig} 
\usepackage{mathtools}
\usepackage{tabu}
\usepackage{amsthm}
\usepackage{empheq}
\usepackage{amsopn}
\usepackage{dsfont}
\usepackage{wrapfig}
\usepackage[font=small]{caption}
\usepackage{units}
\usepackage[symbol]{footmisc}
\usepackage{todonotes}
\usepackage{ wasysym }
\usepackage{amssymb}

\usepackage{kvsetkeys}% provides \comma@parse
\usepackage{etexcmds}% provides \etex@unexpanded

\makeatletter
\newcount\arg@count
\newcommand{\arg@parser}[1]{%
  \advance\arg@count\@ne
  \expandafter\let\csname arg\romannumeral\arg@count\endcsname\comma@entry
}
\newcommand\res[1]{% Default is empty and will be configured later
  \arg@count=\z@
  \comma@parse{ \lambda,A }\arg@parser % Default values
  \arg@count=\z@
  \comma@parse{#1}\arg@parser
  \ifnum\arg@count>2 %
    \@latex@error{Too many arguments}{%
      The macro \string\mycmd\space got \the\arg@count\space
       arguments,\MessageBreak
      but expected are 2 arguments.\MessageBreak
      \@ehd
    }%
  \fi
  \edef\process@me{%
    \noexpand\@res
    {\etex@unexpanded\expandafter{\argi}}%
    {\etex@unexpanded\expandafter{\argii}}%
  }%
  \process@me
}
\newcommand{\@res}[2]{%
  \ensuremath\left( #1 - #2 \right)^{-1}
}
\makeatother

\makeatletter
\newcount\arg@count
\newcommand\p[1]{% Default is empty and will be configured later
  \arg@count=\z@
  \comma@parse{  }\arg@parser % Default values
  \arg@count=\z@
  \comma@parse{#1}\arg@parser
  \ifnum\arg@count=2 %
  \else
    \@latex@error{Wrong number of mandatory arguments}{%
      The macro \string\p\space got \the\numexpr\arg@count-2\relax\space
      mandatory arguments,\MessageBreak
      but expected are 3 mandatory arguments.\MessageBreak
      \@ehd
    }%
  \fi
  \edef\process@me{%
    \noexpand\@p
    {\etex@unexpanded\expandafter{\argi}}%
    {\etex@unexpanded\expandafter{\argii}}%
  }%
  \process@me
}
\newcommand{\@p}[2]{%
  \ensuremath \left\langle #1 , #2 \right\rangle
}
\makeatother

\allowdisplaybreaks

\theoremstyle{definition}
\newtheorem{de}{Definition}[section]

\theoremstyle{plain}
\newtheorem{prop}[de]{Proposition}
\newtheorem{lemma}[de]{Lemma}
\newtheorem{theorem}[de]{Theorem}
\newtheorem{corollary}[de]{Corollary}

\numberwithin{equation}{section}
\theoremstyle{remark}
\newtheorem{remark}[de]{Remark}

\definecolor{dkred}{rgb}{0.7,0,0}
\definecolor{dkgreen}{rgb}{0,0.4,0}

\newcommand{\frank}[1]{{ #1}}

\title{\sc Norm-Resolvent Convergence in Perforated Domains}
\author{
K. Cherednichenko\thanks{\frank{
	Department of Mathematical Sciences\newline
	University of Bath\newline
	Claverton Down\newline
	Bath, BA2 7AY\newline
	United Kingdom\newline
	Email: K.Cherednichenko@bath.ac.uk }
}
, P. Dondl\thanks{\frank{
	Abteilung f\"ur angewandte Mathematik\newline 
	Albert-Ludwigs-Universit\"at Freiburg \newline
	Hermann-Herder-Str. 10\newline
	79104 Freiburg i. Br \newline
	Germany \newline
	Email: \{patrick.dondl, frank.roesler\}@mathematik.uni-freiburg.de}
}
, F. R\"osler\footnotemark[2]
}
\begin{document}

\def\l{\mathopen}
\def\r{\mathclose}

\maketitle  

\begin{abstract}
	For several different boundary conditions (Dirichlet, Neumann, Robin), we prove norm-resolvent convergence for the operator $-\Delta$ in the perforated domain $\Om\setminus \bigcup_{ i\in 2\eps\mathbb Z^d }B_{a_\eps}(i),$ $a_\varepsilon\ll\varepsilon,$ to the limit operator $-\Delta+\mu_{\iota}$ on $L^2(\Om)$, where $\mu_\iota\in\mathbb C$ is a constant depending on the choice of boundary conditions. 
	%We consider Dirichlet, Neumann and Robin boundary conditions.
		
	This is an improvement of previous results \cite{CM}, \cite{Kaizu}, which show \emph{strong} resolvent convergence. In particular, our result implies Hausdorff convergence of the spectrum of the resolvent for the perforated domain problem.
\end{abstract}

\section{Introduction}\label{intro}

In this article we study the following homogenisation problems labelled by $\iota\in\{{\rm D}, {\rm N}, \alpha\}$ (``D" for Dirichlet, ``N'' for Neumann, and ``$\alpha$" for Robin).
\frank{Let $\Om\subset\R^d,$ $d\ge 2,$ be open (bounded or unbounded) with $C^2$ boundary. For unbounded domains $\Omega$ we assume translation invariance, i.e.,  $\Omega + z = \Omega$ for any $z\in \mathbb{Z}^d$.} Let $\alpha\in\mathbb C\setminus\{0\},\;\operatorname{Re}(\alpha)\geq 0$ and denote $\Om_\eps :=\Om\setminus \bigcup_{ i\in L_\eps }B_{a_\eps}(i)$ where $\varepsilon\in(0,1)$, $B_{a_\eps}(i)$ is the ball of radius 
\begin{equation}\label{radius}
	a_\eps^{\rm D}=\begin{cases}
		\eps^{\nicefrac{d}{(d-2)}}, &d\ge 3,\\[0.3em]
	{\rm e}^{-1/\eps^2}, &d=2,
	\end{cases}\quad\quad\quad
	a_\eps^{\rm N}= o(\eps)\ \ (\varepsilon\to0),\qquad\quad a_\eps^\alpha = \eps^{\nicefrac
	{d}{(d-1)}}.
\end{equation}
	centered at the point $i\in L_\eps,$ and 
	\begin{equation}
	L_\eps :=\{i\in 2\eps\mathbb Z^d : \operatorname{dist}(i,\del\Om)>\eps\}.
	\label{Ls}
	\end{equation}
	 Consider the boundary value problems 
\begin{equation}\label{Dirichlet_problem}\tag{Dir}
\begin{cases}
	 (-\Delta+1)u^\eps &= f \text{ in }\Om_{\eps},\\[0.3em]
	\hfill u^\eps &= 0 \text{ on } \del\Om_\eps,
\end{cases}
\end{equation} 

\begin{equation}\label{Neumann_problem}\tag{Neu}
\begin{cases}
	 (-\Delta+1)u^\eps &= f \text{ in }\Om_{\eps},\\[0.3em]
	\hfill \del_\nu u^\eps &= 0 \text{ on } \del\Om_\eps,
\end{cases}
\end{equation}

\begin{equation}\label{Robin_problem}\tag{Rob}
\begin{cases}
	 (-\Delta+1)u^\eps &= f \text{ in }\Om_{\eps},\\[0.3em]
	\hfill \del_\nu u^\eps + \alpha u &= 0 \text{ on } \del\Om_\eps,
\end{cases}
\end{equation}
{\it i.e.} the resolvent problem for the Laplacian, subject to the Dirichlet, Neumann and Robin boundary conditions, respectively.
It is easy to see, using the Lax-Milgram theorem, that for all $\eps\in(0,1)$ each of these problems has a unique weak solution $u^\eps$. It is  a classical question, which we refer to as the homogenisation problem, whether the family of solutions to \eqref{Dirichlet_problem}, \eqref{Neumann_problem}, \eqref{Robin_problem}, obtained by varying the parameter $\varepsilon,$ converges in the sense of the $L^2$-norm to a function $u\in L^2(\Om)$ as $\eps\to 0$ and whether the limit function $u$ solves, in a reasonable sense, some PDE whose form is independent of the right-hand side datum $f.$ 
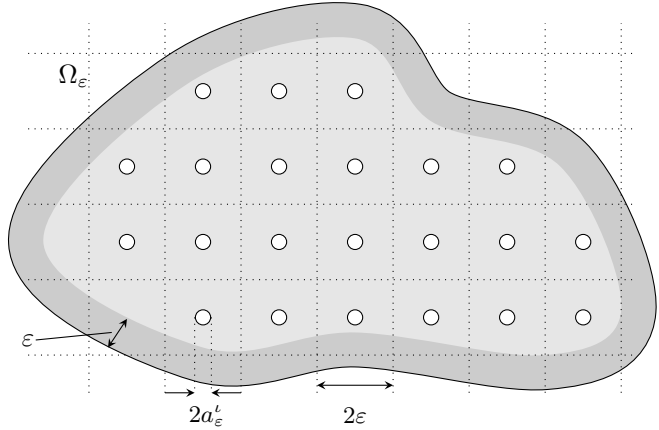
\begin{figure}[htbp]
	\begin{center}
	%\documentclass[11pt,a4paper]{article} 
%
%\usepackage{etex}
%\usepackage[english]{babel}
%\usepackage{tikz}
%\usepackage{pgfplots}
%\usetikzlibrary{fadings,fpu,arrows,decorations.pathmorphing,backgrounds,positioning,fit,petri,patterns}
%\usepackage[utf8]{inputenc}
%\usepackage{graphicx}
%\usepackage[colorlinks=true,linkcolor=black,citecolor=black]{hyperref}
%\usepackage{amsfonts}
%\usepackage{enumerate}
%\usepackage{booktabs}
%\usepackage{array} 
%\usepackage{paralist} 
%\usepackage{subfig} 
%\usepackage{mathtools}
%\usepackage{tabu}
%\usepackage{amsthm}
%\usepackage{empheq}
%\usepackage{amsopn}
%\usepackage{dsfont}
%\usepackage{wrapfig}
%\usepackage[font=small]{caption}
%
%\begin{document}

\begin{tikzpicture}[>=stealth]

	\filldraw[double distance=9mm,fill = gray!20, draw=black, double=gray!40] plot [smooth cycle, tension=0.6] coordinates {(6.5,-1.5) (4,-1.2) (2,-1.4) (-0.1,0) (1.8,2.1) (4,2.7) (5,1.6) (6.7,1) (7.5,-0.9)};
	\fill[fill = gray!20] plot [smooth cycle, tension=0.6] coordinates 
		{(6.5,-1.5) (4,-1.2) (2,-1.4) (-0.1,0) (1.8,2.1) (4,2.7) (5,1.6) (6.7,1) (7.5,-0.9)};

	\draw (0.3,2.2)node[]{${\Omega_\varepsilon}$} ;
	
	\filldraw[fill = white] (1,0) circle (1mm) ;
	\filldraw[fill = white] (2,0) circle (1mm) ;
	\filldraw[fill = white] (3,0) circle (1mm) ;
	\filldraw[fill = white] (4,0) circle (1mm) ;
	\filldraw[fill = white] (5,0) circle (1mm) ;
	\filldraw[fill = white] (6,0) circle (1mm) ;
	\filldraw[fill = white] (7,0) circle (1mm) ;

	\filldraw[fill = white] (1,1) circle (1mm) ;
	\filldraw[fill = white] (2,1) circle (1mm) ;
	\filldraw[fill = white] (3,1) circle (1mm) ;
	\filldraw[fill = white] (4,1) circle (1mm) ;
	\filldraw[fill = white] (5,1) circle (1mm) ;
	\filldraw[fill = white] (6,1) circle (1mm) ;
	
	\filldraw[fill = white] (3,2) circle (1mm) ;
	\filldraw[fill = white] (4,2) circle (1mm) ;
	\filldraw[fill = white] (2,2) circle (1mm) ;

	\filldraw[fill = white] (2,-1) circle (1mm) ;
	\filldraw[fill = white] (3,-1) circle (1mm) ;
	\filldraw[fill = white] (4,-1) circle (1mm) ;
	\filldraw[fill = white] (5,-1) circle (1mm) ;
	\filldraw[fill = white] (6,-1) circle (1mm) ;
	\filldraw[fill = white] (7,-1) circle (1mm) ;

	\draw[dotted, shift={(0.5,0.5)}] (-0.8,-2.5) grid (7.2,2.4);
	
	\draw[<->] (0.75,-1.4) -- (1,-1.02);
	\draw (0.8,-1.2) -- (-0.2,-1.3);
	\draw (-0.3,-1.3) node[]{{$\varepsilon$}};
	
	\draw[->] (1.5,-2) -- (1.89,-2);
	\draw[->] (2.5,-2) -- (2.11,-2);
	\draw[dotted] (1.89,-1) -- (1.89,-2);
	\draw[dotted] (2.11,-1) -- (2.11,-2);
	\draw (2.05,-2.3) node[]{\small${2a_\varepsilon^\iota}$};
	
	\draw[<->,shift={(-1,0)}] (4.5,-1.9) -- (5.5,-1.9);
	\draw[shift={(-1,0)}] (5,-2.3) node[]{\small${2\varepsilon}$};
	
\end{tikzpicture}

%\end{document}
	\end{center}
	\caption{ Sketch of the perforated domain.}
\end{figure}

Homogenisation problems of this type have been studied extensively for a long time \cite{CM,RT,MK,Kaizu}. For example, results by Cioranescu \& Murat and Kaizu give a positive answer to the previous question for all three choices of boundary conditions at least in the case of \emph{bounded} domains. In fact, they showed that the solutions of \eqref{Dirichlet_problem}, \eqref{Robin_problem}, \eqref{Neumann_problem} converge strongly in $L^2(\Om)$ to the solution $u\in H^1(\Om)$ of $(-\Delta+1+\mu_\iota)u=f$, where 
	\begin{equation}\label{mu_iota}
		\mu_\iota = \begin{cases}
 			\dfrac{\pi}{2}, &\iota = D,\ d=2,\\[0.8em]
 			\dfrac{(d-2)S_d}{2^d}, &\iota = D,\ d\geq 3,\\[0.8em]
 			0, &\iota = N,\\[0.8em]
			\dfrac{\alpha S_d}{2^d}, &\iota = \alpha.
 		\end{cases}
	\end{equation}
	where $S_d$ denotes the surface area of the unit ball in $\R^d$.
	
In this article we attempt to improve this result in two directions. First, we show the above convergence not only in the strong sense, but in the \emph{norm-resolvent sense} (that is, the right-hand side $f$ is allowed to depend on $\eps$). Second, our result is then extended to the case of unbounded domains. As a corollary, we obtain a statement about the convergence of the spectra of the perforated domain problems \eqref{Dirichlet_problem}, \eqref{Neumann_problem}, \eqref{Robin_problem} as $\varepsilon\to0.$

The paper is organised as follows. In Section \ref{setting} we will briefly give a more precise formulation of the problem and include previous results. In Section \ref{results} we will state our main result and its implications. Sections \ref{uniformity}, \ref{exponentialdecay} and \ref{decomposition} contain the proof of the main theorem and in Section \ref{semigroup} we consider implications of our main theorem for the semigroup generated by the Robin Laplacian. Section \ref{conclusion} contains a brief conclusion and discusses open problems.

\section{Geometric setting and previous results}\label{setting}

\noindent As above, assume $d\geq 2,$ and let
 $$
 T_\eps\;:=\;\bigcup_{ i\in L_\eps}T_i^\eps,\quad\quad T_i^\eps:=B_{a_\eps^\iota}(i),\ \  i\in L_\eps,
 %\;:=\bigcup_{ i\in L_\eps^\iota}B_{a_\eps^\iota}(i),
 $$ 
where $a_\eps^\iota,$ $L_\eps$ as in \eqref{radius}, \eqref{Ls}. Denote $\Om_\eps:=\Om\setminus T_\eps$. We also denote $B_i^\eps:=B_\eps(i)$ and $P_i^\eps:=\eps[-1,1]^d+i$ for $i\in L_\eps.$ Constants independent of $\eps$ will be denoted $C$ and may change from line to line. Note that our assumptions on $\Om$ ensure that the set $\{\phi|_\Om : \phi\in C_0^\infty(\R^d)\}$ is dense in $H^1(\Om)$ (cf. \cite[Cor. 9.8]{brezis}) in the cases $\iota = N,\alpha$.

Moreover, since we are dealing with varying spaces $L^2(\Om_\eps)$, it is convenient to define the identification operators
\begin{align}
	J_\eps &: L^2(\Om_\eps) \to L^2(\Om),
	&J_\eps f(x) &= \begin{cases}
		f(x), & x\in\Om_\eps, \\[0.1em]
		0, & x\in \Om\setminus\Om_\eps
	\end{cases}\label{eq:Jdef}
\\[1mm]
	I_\eps &: \frank{L^2(\Om) \to L^2(\Om_\eps)},
	&I_\eps g(x) &= g|_{\Om_\eps}\label{eq:Idef}
\\[1mm]
	\mathcal T_\eps &: H^1(\Om_\eps) \to H^1(\Om),
	&\mathcal T_\eps u &= \begin{cases}
		u & \text{ in }\Om_\eps, \\[0.1em]
		v & \text{ in } T_\eps,
	\end{cases}\label{Tdef}
\end{align}
where $v$ is the harmonic extension of $u$ into the holes, i.e. 
\begin{equation}
	\begin{cases}
		\Delta v &= 0\, \text{ in }\,T_\eps, \\[0.1em]
		\phantom{\Delta} v &= u\, \text{ on }\, \del T_\eps.
	\end{cases}
\end{equation}
\begin{lemma}
	For $I_\eps,J_\eps$ as above, one has
	\begin{align}
		I_\eps J_\eps &= \mathrm{id}_{L^2(\Om_\eps)}\\
		\|J_\eps I_\eps - \mathrm{id}_{L^2(\Om)}\|_{\mathcal L(H^1(\Om),L^2(\Om))} &\to 0. \label{eq:JI-1}
	\end{align} 
	Moreover, $\|I_\eps\|_{\mathcal L(L^2(\Om_\eps),L^2(\Om))},\, \|J_\eps\|_{\mathcal L(L^2(\Om),L^2(\Om_\eps))}$ are uniformly bounded in $\eps$.
\end{lemma}
\begin{proof}
	The only nontrivial statement is \eqref{eq:JI-1}. To prove this, let $f\in H^1(\Om)$. Then $\|f-J_\eps I_\eps f\|_{L^2(\Om)} = \|f\|_{L^2(T_\eps)}$. To show that this quantity converges to 0 uniformly in $f$, denote $Q_k:=[0,1)^d+k$ for $k\in\mathbb Z^d$ a cube shifted by $k$, so that $\R^d = \bigcup_{k\in\mathbb Z^d}Q_k$. Then we have
	\begin{align*}
		\|f\|_{L^2(T_\eps)}^2 &= \sum_{k\in\mathbb Z^d} \|f\|_{L^2(Q_k\cap T_\eps)}^2\\
		&\leq \sum_{k\in\mathbb Z^d} \left\|1\right\|_{L^{2p}(Q_k\cap T_\eps)}^2 \left\|f\right\|_{L^{2q}(Q_k\cap T_\eps)}^2
	\end{align*}
	for $p,q>1$ with $p^{-1}+q^{-1}=1$, by H\"older's inequality. Since $f\in H^1(\Om)$, we can use the Gagliardo-Sobolev-Nierenberg inequality to conclude (for suitable $q$) that
	\begin{align*}
		\|f\|_{L^2(T_\eps)}^2 &\leq \left\|1\right\|_{L^{2p}(Q_0\cap T_\eps)}^2 \sum_{k\in\mathbb Z^d}  \left\|f\right\|_{L^{2q}(Q_k\cap T_\eps)}^2\\
		&\leq \left\|1\right\|_{L^{2p}(Q_0\cap T_\eps)}^2 \sum_{k\in\mathbb Z^d}  C \left\|f\right\|_{H^1(Q_k)}^2\\
		&= |Q_0\cap T_\eps|^{\nicefrac 1 p}  C \left\|f\right\|_{H^1(\Om)}^2
	\end{align*}
	with some suitable $p>0$. Since $|Q_0\cap T_\eps|\to 0$ as $\eps\to 0$ (cf. the definition of $a_\eps^\iota$, \eqref{radius}), the desired convergence follows.
\end{proof}
\begin{lemma}\label{harmext}
	The harmonic extension operator $\mathcal T_\eps$ satisfies
	\begin{enumerate}[(i)]
		\item 
		$\limsup_{\eps\to 0} \|\mathcal T_\eps\|_{H^1(\Om_\eps)\to H^1(\Om)}<\infty.$
		\item
		There exists $C>0$ such that $\|\mathcal T_\eps w\|_{H^1(P_i^\eps)} \leq C\|w\|_{H^1(P_i^\eps)}$ for all $w\in H^1(\Om_\eps)$ and $i\in L_\eps$.
		\item 
		For any sequence $w_\eps$ such that $\limsup_{\eps\to 0}\|w_\eps\|_{H^1(\Om_\eps)}<\infty$ one has $ \|\mathcal T_\eps w_\eps -  J_\eps w_\eps\|_{L^2(\Om)} \to 0.$
	\end{enumerate}
\end{lemma}
\begin{proof}
	See \cite{Kaizu}, \cite[p.\,40]{RT}.
\end{proof}

In the above geometric setting, we will study the linear operators $A_\eps^\iota,$ $\iota={\rm D}, {\rm N}, \alpha$ in \frank{$L^2(\Omega_\eps),$} defined by the differential expression
%\begin{align*}
%	A_\eps^\iota :&= 
	$-\Delta+1,$
	%\text{ on } 
	%L^2(\Om_\eps),$
%\end{align*}
with (dense) domains
\begin{align*}
	\mathcal D(A_\eps^{\rm D}) &= H^1_0(\Om_\eps)\cap H^2(\Om_\eps),\\[0.3em]
	\mathcal D(A_\eps^{\rm N}) &= \bigl\{u\in H^2(\Om_\eps):\del_\nu u=0\text{ on }\del\Om_\eps\bigr\},
	\\[0.3em]
	\mathcal D(A_\eps^\alpha) &=\bigl\{u\in H^2(\Om_\eps):\del_\nu u+\alpha u = 0 \text{ on }\del\Om_\eps\bigr\},
\end{align*}
respectively, and the linear operators $A^\iota$ in $L^2(\Omega_\eps)$ defined by the expression
%\begin{align*}
	%A^\iota :&= 
	$-\Delta+1+\mu_\iota,$
	%\quad\text{ on } L^2(\Om)
%\end{align*}
with domains
\begin{align*}
	\mathcal D(A^{\rm D}) &= H^1_0(\Om)\cap H^2(\Om),\\[0.3em]
	\mathcal D(A^{\rm N}) &= \bigl\{u\in H^2(\Om):\del_\nu u=0\text{ on }\del\Om\bigr\},\\[0.3em]
	\mathcal D(A^\alpha) &=\bigl\{u\in H^2(\Om):\del_\nu u+\alpha u = 0 \text{ on }\del\Om\bigr\},
\end{align*}
respectively, where $\mu_\iota,$ $\iota={\rm D}, {\rm N}, \alpha,$  are defined in \eqref{mu_iota}.\begin{remark}
In the case when $d\ge 3$ one has the characterisation 
\begin{equation}
	\mu_{\rm D} = \dfrac{1}{2^d}\inf\biggl\{\int_{\R^d\setminus B_1(0)}\vert\nabla u\vert^2,\ \ \ u\in H^1({\mathbb R}^d),\ u=1{\rm \ on\ }B_1(0)\biggr\}.
	\label{capacity_form}
\end{equation} 
Note that the factor $1/2^{d}$ arises from the fact that the unit cell is of size $2\eps$.	
\end{remark}
Using the notation above, we recall the following classical results.
\begin{theorem}[\cite{CM}]\label{CMresult}
	Let $\Om\subset\R^d$ be open (bounded or unbounded). Suppose that $f\in L^2(\Om),$ and let $u^\eps$ and $\tilde u$ be the solutions to
	\begin{align*}
		(-\Delta+1)u^\eps &= f ,\qquad u^\eps\in H^1_0(\Om_{\eps}),\\[0.25em]
		(-\Delta+1+\mu_{\rm D})\tilde u &= f, \qquad \tilde u\in H^1_0(\Om).
	\end{align*}
	Then  
	$
		J_\eps u^\eps\xrightharpoonup{\eps\to 0}\tilde u$
		%\quad\text{ in }\;
	in	
	$
	H^1_0(\Om).
	$
\end{theorem}
\begin{theorem}[\cite{Kaizu}]\label{Kaizuresult}
	Let $\Om\subset\R^d$ be open (bounded or unbounded), and suppose that $\del\Omega$ is smooth. Suppose also that $f\in L^2(\Om),$ and let $u^\eps$ and $\tilde u$ be the solutions to
	\begin{align*}
		(-\Delta+1)u^\eps &= f ,\qquad u^\eps\in \mathcal D(A_\eps^{\alpha, {\rm N}}),\\[0.25em]
		(-\Delta+1+\mu_{\alpha, {\rm N}})\tilde u &= f, \qquad \tilde u\in \mathcal D(A^{\alpha, {\rm N}}).
	\end{align*}
	Then one has 
	$$
		\frank{\mathcal T_\eps} u^\eps\xrightharpoonup{\eps\to 0}\tilde u\quad\text{ in }\;H^1(\Om).
	$$
\end{theorem}
\begin{proof}[Proof of Theorems \ref{CMresult} and \ref{Kaizuresult}]
	The results are obtained by following the proofs of \cite[Thm 2.2]{CM}, \cite[Thm 2]{Kaizu}. Note that the weak convergence in $H^1(\Om)$ is immediately obtained also for unbounded domains (and complex $\alpha$).
\end{proof}
An important ingredient in the proofs are auxiliary functions $w_\epsilon^{\iota}\in W^{1,\infty}(\R^d)$ defined, for each $\varepsilon\in(0,1),$ as the solution to the problems
\begin{equation}\label{wepsilon}
	w_\eps^{\rm N} \equiv 1, \qquad\quad
	\begin{cases}
		\phantom{\Delta}w_\eps^{\rm D} = 0 & \text{ in } T_i^\eps,\\[0.2em]
		\Delta w_\eps^{\rm D} = 0 &\text{ in } B_i^\eps\setminus T_i^\eps,\\[0.2em]
		\phantom{\Delta}w_\eps^{\rm D} = 1 &\text{ in }P_i^\eps\setminus B_i^\eps,\\[0.2em]
		\phantom{\Delta}w_\eps^{\rm D} &\text{continuous,}
	\end{cases} \qquad\quad
	\begin{cases}
		\del_\nu w_\eps^{\alpha} + \alpha w_\eps^{\alpha} = 0 & \text{ on } \del T_i^\eps,\\[0.2em]
		\Delta w_\eps^\alpha = 0 &\text{ in } B_i^\eps\setminus T_i^\eps,\\[0.2em]
		\phantom{\Delta}w_\eps^\alpha = 1 &\text{ in }P_i^\eps\setminus B_i^\eps,\\[0.2em]
		\phantom{\Delta}w_\eps^\alpha &\text{continuous,}
	\end{cases}
\end{equation}
used as a test function in the weak formulation of the problems 
\eqref{Dirichlet_problem}, \eqref{Neumann_problem}, \eqref{Robin_problem}. 
	These functions were used in \cite{CM,Kaizu} as test functions to prove strong convergence of solutions. They are ``optimal'' in the sense that they minimise the energy in annular regions around the holes. In the Dirichlet case, the function $w_\eps^{\mathrm{D}}$ is nothing but the potential for the capacity $\operatorname{cap}\bigl(B_\eps(i);B_{a_\eps^{\mathrm D}}(i)\bigr)$. It can be shown that one has the convergences
	\begin{align}
	\begin{rcases}
		\phantom{\mathcal T_\eps }w_\eps^{\mathrm D} &\rightharpoonup 1\\
		\mathcal T_\eps w_\eps^{\mathrm \alpha}&\rightharpoonup 1
	\end{rcases}
	&\text{ weakly in } H^{1}(\Om)\\[1mm]
	\begin{rcases}
		\qquad\qquad\qquad\qquad\quad-\Delta w_\eps^{\mathrm D} &\to\; \mu_{\mathrm D}\; \\
		-\nabla\cdot (\chi_{\Om_\eps}\nabla w_\eps^\alpha) + \alpha w_\eps^\alpha \delta_{\del T_\eps} &\to\; \mu_\alpha \;
	\end{rcases}
	&\text{ strongly in } H^{-1}(\Om)\label{eq:weps_convergence}
	\end{align}
	as $\eps\to 0$, where $\delta_{\del T_\eps}$ denotes the Dirac measure on the boundary of the holes (for a proof of the above facts, see \cite[Lemma 2.3]{CM} and \cite[Section 3]{Kaizu}).

\section{Main results}\label{results}

In what follows we prove the following claim.
% we are going to prove 

\begin{theorem}\label{mainth}
	Let $J_\eps,A_\eps^\iota,A^\iota$ be defined as in the previous section. Then for $\iota\in\{{\rm D}, {\rm N}, \alpha\}$ one has 
	\begin{align}\label{eq:maineq}
		\bigl\|J_\eps (A_\eps^\iota)^{-1}-(A^\iota)^{-1}J_\eps\bigr\|_{\mathcal L(L^2(\Om_\eps),\, L^2(\Om))}\to 0\qquad (\eps\to 0),	
	\end{align}

	that is, the operator sequence $A_\eps^\iota$ converges to $A^\iota$ in the norm-resolvent sense.
\end{theorem}

%\begin{corollary}
%	If $A_\eps, A$ are as in Theorem \ref{mainth}, then 
%	\begin{align}\label{eq:AI-IA}
%		\bigl\|(A_\eps^\iota)^{-1}I_\eps - I_\eps(A^\iota)^{-1}\bigr\|_{\mathcal L(L^2(\Om),	L^2(\Om_\eps))}\to 0,
%	\end{align}
%	where $I_\eps$ is as in \eqref{eq:Idef}.
%\end{corollary}
%\begin{proof}
%	For notational convenience, denote $R_\eps:=(A_\eps^\iota)^{-1}$ and $R:=(A^\iota)^{-1}$. A quick calculation shows that 
%	\begin{align*}
%		R_\eps I_\eps - I_\eps R &= I_\eps(J_\eps R_\eps - R J_\eps)I_\eps - (I_\eps J_\eps - 1)R_\eps I_\eps\\
%		&= I_\eps(J_\eps R_\eps - R J_\eps)I_\eps,
%	\end{align*}
%	since $I_\eps J_\eps = \mathrm{id}_{L^2(\Om_\eps)}$. Hence
%	\begin{align*}
%		\|R_\eps I_\eps - I_\eps R\|_{\mathcal L(L^2(\Om),L^2(\Om_\eps))}
%		&\leq \left\|I_\eps\right\|_{\mathcal L(L^2(\Om),L^2(\Om_\eps)}^2 \left\|J_\eps R_\eps - R J_\eps\right\|_{\mathcal L(L^2(\Om_\eps),L^2(\Om))}\\
%		&\to 0
%	\end{align*}
%	as $\eps\to 0$, by \eqref{eq:maineq} and uniform boundedness of $\|I_\eps\|_{\mathcal L(L^2(\Om_\eps),L^2(\Om))}$.
%\end{proof}

We note an important consequence of the above theorem.

\begin{corollary}\label{SpCon}
	If the boundary conditions are selfadjoint, then for all bounded open $V\subset\mathbb C,$ one has $\sigma(A_\eps^\iota)\cap V\xrightarrow{\eps\to 0}\sigma(A^\iota)\cap V$ in the Hausdorff sense.\footnote{For the definition of Hausdorff convergence, see {\it e.g.} \cite{Rockafellar}.}
\end{corollary}
\begin{proof}
	First, note that the spectra of $A_\eps^\iota$ converge to that of $A^\iota,$ in the sense that for each compact $K\subset \rho(A^\iota)$ there exists $\eps_0>0$ such that $K\subset\rho(A_\eps^\iota)$ for all $\eps\in(0,\eps_0)$. The proof of this is obtained by combining the proofs of Lemma 3.11, Theorem 3.12 and Corollary 3.14 in \cite{MNP}.
	On the other hand, the implication $K\subset\rho(A_\eps)\,\forall\eps>0\Rightarrow K\subset\rho(A)$ is standard for selfadjoint operators (cf. \cite[Th. VIII.24]{RS}).
	
	Together these two facts imply the assertion. Indeed, for given $V\subset\C$, let $\delta>0$ and define the compact set $K:=\overline{V}\setminus U_\delta(\sigma(A^\iota))$, where $U_\delta(\cdot)$ denotes the open $\delta$-neighbourhood. By the above we have that $K\subset \rho(A_\eps^\iota)$ for $\eps$ small enough. This shows that $V\cap\sigma(A_\eps^\iota)\subset V\cap U_\delta(\sigma(A^\iota))$. 
	
	To see the converse inclusion $V\cap \sigma(A^\iota)\subset V\cap U_\delta(\sigma(A^\iota_\eps))$, let us argue by contradiction and suppose that there exists $\delta_0>0$ such that
	\begin{equation*}
		V\cap\sigma(A^\iota)\nsubseteq V\cap U_{\delta_0}(\sigma(A^\iota_\eps))\qquad\forall\eps>0.
	\end{equation*}
	By this assumption, there exists a sequence $(\lambda_\eps)\subset  V\cap\sigma(A^\iota)$ such that  $\dist(\lambda_\eps,\sigma(A^\iota_\eps))\ge \delta_0$ for all $\eps>0$. Since $(\lambda_\eps)$ is bounded, we can extract a subsequence $\lambda_{\eps'}\to\lambda_0\in \overline V\cap\sigma(A^\iota)$.
	It follows that
	\begin{equation*}
		\overline{B_{\nicefrac{\delta_0}{2}}(\lambda_0)}\subset\rho(A^\iota_{\eps'})\qquad \text{for almost all }\eps'>0.
	\end{equation*}
	By norm-resolvent convergence of $(A_{\eps'}^\iota)$, we conclude that $\overline{B_{\nicefrac{\delta_0}{2}}(\lambda_0)}\subset\rho(A^\iota)$, contradicting the fact that $\lambda_0\in\sigma(A^\iota)$.
	
	Since $\delta>0$ was arbitrary, the desired Hausdorff convergence follows.
\end{proof}
In particular, this corollary shows that (if $\re(\mu_\iota)>0$) a spectral gap opens for $A^\iota_\eps$ between $0$ and $\re(\mu_\iota)$.

\begin{remark}
	We note that our assumption on the spherical shape of the holes was made for the sake of definiteness, however, our results easily generalise to more general geometries as detailed in \cite[Th. 2.7]{CM}. Moreover, our results are also valid  for more general elliptic operators $\operatorname{div}(A\nabla)$ with continuous coefficients $A$ (cf. \cite{CM}).
\end{remark}

\section{Uniformity with respect to the right-hand side}\label{uniformity}

In this section we prove that the result of Theorems \ref{CMresult}, \ref{Kaizuresult} hold in a strengthened form, namely, uniformly with respect to the right-hand side $f$. More precisely, the following  holds.
\begin{theorem}\label{uniright}
	Suppose that $\eps_n\searrow 0,$ 
	$f_n\in L^2(\Om_{\eps_n}),$ $n\in{\mathbb N},$ with $\|f_n\|_{L^2(\Omega_{\eps_n})}\leq 1,$ and let $u_n^\iota$ and $\tilde u_n^\iota$ be the solutions to the problems ($\iota\in\{{\rm D}, {\rm N}, \alpha\}$)
	\begin{align}
		(-\Delta+1)u_n^\iota &= f_n ,\qquad u_n^\iota\in \mathcal D(A_{\eps_n}^\iota),\label{uniright1}\\[0.25em]
		(-\Delta+1+\mu_\iota)\tilde u_n^\iota &= J_{\eps_n} f_n, \qquad \tilde u_n^\iota\in \mathcal D(A^\iota).\label{uniright2}
	\end{align}
	Then for every bounded, open $K\subset\Omega$ one has 
	\frank{
	\begin{align*}
		J_{\eps_n}u_n^\iota-\tilde u_n^\iota &\rightarrow 0 \qquad\text{strongly in }\;L^2(K),\\[0.25em]
			J_{\eps_n}\nabla u^\iota_n - \nabla \tilde u^\iota_n &\rightharpoonup 0 \qquad \text{weakly in }L^2(K),
	\end{align*}
	}
	for $\iota\in \{{\rm D}, {\rm N},\alpha \}$.
\end{theorem}

\begin{proof}
	We have the following {\it a priori} estimates (note Lemma \ref{harmext}):
	\begin{align*}
		\|\mathcal{ T}_{\eps_n}u_n^{\alpha, {\rm N}}\|_{H^1(\Omega)} &\leq C\|J_{\eps_n} f_n\|_{L^2(\Omega)},\\[0.3em]
		\|J_{\eps_n} u_n^{\rm D}\|_{H^1(\Omega)} &\leq C\|J_{\eps_n}f_n\|_{L^2(\Omega)}, \\[0.3em]
		\|\tilde u_n^\iota\|_{H^1(\Omega)} &\leq C\|J_{\eps_n}f_n\|_{L^2(\Omega)}\quad\forall\iota\in\{{\rm D}, {\rm N}, \alpha\}.
	\end{align*}
	Thus, there exists a subsequence (still indexed by $n$) and $u^\iota,\tilde u^\iota\in H^1(\Om)$ such that 
	\begin{align}
		\begin{rcases}
			J_{\eps_n}u_n^{\rm D} &\xrightharpoonup{n\to\infty}u^{\rm D} \; \\[0.25em]
			\mathcal T_{\eps_n}u_n^{\alpha, {\rm N}}&\xrightharpoonup{n\to\infty}u^{\alpha, {\rm N}} \; \\[0.25em]
			\qquad \tilde u_{n}^\iota&\xrightharpoonup{k\to\infty}\tilde u^\iota,\ \ \ \iota\in\{{\rm D}, {\rm N}, \alpha\} \;
		\end{rcases}
		\text{ weakly in } H^1(\Om).
		\label{two_convergences}
	\end{align}
	Note that that for every bounded $K\subset\Omega$ the convergence statements (\ref{two_convergences}) are strong in $L^2(K)$. In particular, employing Lemma \ref{harmext} (i), (iii) we immediately obtain
	\frank{
	\begin{align}
			J_{\eps_n}u_n^\iota &\to u^\iota \qquad\;\text{ strongly in }L^2(K),\nonumber \\[0.2em]
			J_{\eps_n}\nabla u^\iota_n &\rightharpoonup \nabla u^\iota \qquad \text{weakly in }L^2(K).\nonumber \\[0.4em]
			u_n^\iota &\to \tilde u^\iota \qquad\;\text{ strongly in }L^2(K), \label{strL2} \\[0.2em]
			\nabla u^\iota_n &\rightharpoonup \nabla \tilde u^\iota \qquad \text{weakly in }L^2(K). \label{weakL2}
	\end{align}
	}
	for all $\iota\in\{{\rm D}, {\rm N}, \alpha\}$. Next, choose a further subsequence (still indexed by $n$) such that also $J_{\eps_n}f_n\xrightharpoonup{n\to\infty}f$ weakly in $L^2(\Om),$ where the limit $f$ may depend on the choice of subsequence. 
	\paragraph{\it Dirichlet and Neumann case.} We restrict ourselves to the Dirichlet and Neumann problems first and comment on the Robin problem at the end of the proof. Consider the weak formulations of the problem \eqref{uniright2}, {\it i.e.} 
	\begin{align*}
		\int_{\Om}\overline{\nabla\tilde u_n^\iota}\nabla\phi + (1+\mu_\iota)\int_{\Om}\overline{\tilde u_n^\iota}\phi = \int_{\Om} \overline{f_n}\phi,	
	\end{align*}
	where $\phi\in C_0^\infty(\Om)$ for $\iota={\rm D}$ and $\phi\in C_0^\infty(\R^d)$ for $\iota={\rm N}$. 
	Letting $n\to\infty$ and using the \frank{convergences \eqref{strL2}, \eqref{weakL2}} (with $K=\Om\cap\supp\phi$) we obtain
	\begin{align*}
		\int_{\Om}\overline{\nabla\tilde u^\iota}\,\nabla\phi + (1+\mu_\iota)\int_{\Om}\overline{\tilde u^\iota}\phi = \int_{\Om} \overline{f}\phi.
	\end{align*}
	Next consider the weak formulation of \eqref{uniright1},where we choose the test function $w_{\eps_{n}}^\iota\phi$:
	\begin{align*}
		\int_{\Om_{\eps_n}}\overline{\nabla u_n^\iota}\nabla\big(w_{\eps_n}^\iota\phi\big) + \int_{\Om_{\eps_n}} \overline{u^\iota_n} w_{\eps_n}^\iota\phi = \int_{\Om_{\eps_n}} \overline{f_n} w_{\eps_n}^\iota\phi,
	\end{align*}
	where again $\phi\in C_0^\infty(\Om)$ for $\iota={\rm D}$ and $\phi\in C_0^\infty(\R^d)$ for $\iota= {\rm N}$.
	It follows from the results of \cite{CM,Kaizu} (cf. \eqref{eq:weps_convergence}) that the left and right-hand side of this equation converge to
	$$
		\int_{\Om}\big(\overline{\nabla u^\iota}\nabla\phi + (1+\mu_\iota) \overline{u^\iota}\phi\big)\quad
	%$$
	{\rm and}\quad
	%$$
		\int_{\Om} \overline{f}\phi,
	$$
	respectively. Thus, we obtain 
	$$
		\int_{\Om}\big(\overline{\nabla u^\iota}\,\nabla\phi + (1+\mu_\iota) \overline{u^\iota}\phi\big)=\int_{\Om} \overline{f}\phi,
	$$ 
	and hence $u^\iota$ and $\tilde u^\iota$ are weak solutions to the same equation. Uniqueness of solutions (for all $\iota\in\{{\rm D}, {\rm N}\}$) implies $\tilde u^\iota = u^\iota,$ which shows the assertion for the chosen subsequence.
	
	Finally, applying the above reasoning to every subsequence of $(J_{\eps_n}u_n^\iota-\tilde u_n^\iota)$ yields the result for the whole sequence.
	\frank{
	\paragraph{\it Robin case.} In the Robin case, the above proof remains valid in the interior of $\Om_\eps$, but convergence of the boundary terms 
	$$
		\int_{\del T_{\eps_n}} w^\iota_{\eps_n}u^\iota_n\phi\qquad\text{and}\qquad \int_{\del\Omega} u^\iota_n\phi
	$$
	has to be shown. Convergence of the second term follows since $u^\iota_n\rightharpoonup u^\iota$ in $L^2(\del\Om)$, while convergence of the first term follows from \eqref{eq:weps_convergence}. For details, see \cite{Kaizu}.
	}
\end{proof}

\begin{corollary}
	If the domain $\Om$ is bounded, one has
	$$\bigl\|J_\eps (A^\iota_\eps)^{-1}-(A^\iota)^{-1}J_\eps\bigr\|_{\mathcal L(L^2(\Om_\eps),\, L^2(\Om))}\to 0\qquad (\eps\to 0)$$ 
	for $\iota\in\{{\rm D}, {\rm N}, \alpha\}$, i.e., Theorem \ref{mainth} holds in that case of bounded $\Om$.
\end{corollary}
\begin{proof}
	Since $\Om$ is bounded, the embedding of $H^1(\Om)$ in $L^2(\Om)$ is compact, thus the sequence $J_{\eps_n}u_n^\iota-\tilde u_n^\iota$ from the previous proof has a subsequence converging to 0 strongly in $L^2(\Om)$. Since this can be done for every subsequence of $(J_{\eps_n}u_n^\iota-\tilde u_n^\iota)$, the whole sequence converges to 0. 
	
	Now, choose a sequence $f_n\in L^2(\Om_{\eps_n}),\;\|f_n\|_{L^2(\Om_\eps)}\leq 1$, such that
	$$
		\sup_{\substack{\scriptscriptstyle{f\in L^2(\Om_{\eps_n})} \\ \scriptscriptstyle{\|f\|\leq 1}}} \big\|(J_{\eps_n} (A_\eps^\iota)^{-1}-(A^\iota)^{-1}J_{\eps_n})f \big\|_{L^2(\Om_\eps)}-\f 1 n < \big\|(J_{\eps_n} (A_{\eps_n}^\iota)^{-1}-(A^\iota)^{-1}J_{\eps_n})f_n \big\|_{L^2(\Om_{\eps_n})}.
	$$
	By the above, the right-hand side of this inequality converges to zero, which implies the claim.
\end{proof}
Treating unbounded domains requires further effort. Since we lack compact embeddings in this case, we will have to take advantage of the sufficiently rapid decay of solutions to $(-\Delta+1)u = f$ and a decomposition of the right hand side with a bound on the interactions.

\section{Exponential decay of solutions}\label{exponentialdecay}

We begin with a general result which we assume is classical, but include for the sake of completeness. Let $U\subset\R^d$ open satisfying the strong local Lipschitz condition, $\lambda>\f 1 2$ and consider the problems ({\it cf.} (\ref{Dirichlet_problem}), (\ref{Neumann_problem}), (\ref{Robin_problem}))
\begin{align}\label{GenRob}
	\begin{cases}
	(-\Delta+\lambda)u^{\alpha} &= f\quad\text{ in }U, \\[0.3em]
	\hfill \del_\nu u^\alpha + \alpha u^\alpha &= 0 \quad\text{ on }\del U;
	\end{cases}
\end{align}
\begin{align}\label{GenNeu}
	\begin{cases}
	(-\Delta+\lambda)u^{\rm N} &= f\quad\text{ in }U, \\[0.3em]
	\hfill \del_\nu u^{\rm N} &= 0 \quad\text{ on }\del U;
	\end{cases}
\end{align}
\begin{align}\label{GenDir}
	\begin{cases}
	(-\Delta+\lambda)u^{\rm D} &= f\quad\text{ in }U, \\[0.3em]
	\hfill  u^{\rm D} &= 0 \quad\text{ on }\del U.
	\end{cases}
\end{align}
Let $x_0\in\R^d,$ and define the function $\omega(x)=\cosh(|x-x_0|)$. Then the following statement holds.

\begin{prop}\label{main}
	Let $f\in L^2(U)$, $\supp(f)$ compact. Then each of the problems \eqref{GenRob}--\eqref{GenDir} has a unique weak solution $u^\iota\in H^1(U)$ satisfying 
	\begin{align}
		\int_U |u^\iota|^2\omega\,dx &\leq M\int_U |f|^2\omega\,dx \label{est1}\\
		\int_U |\nabla u^\iota|^2\omega\,dx &\leq M\int_U |f|^2\omega\,dx ,\label{est2}
	\end{align}
	where $M:=\max\left\{ 2,(\lambda-\f 1 2)^{-1} \right\}$.
\end{prop}
We postpone the proof, in order to introduce some notation and prove auxiliary results. First, let us denote $d\mu:=\omega dx$ and introduce the weighted Sobolev spaces $\mathcal H := W^{1,2}(U;\omega)$, $\mathcal H_0 := W^{1,2}_0(U;\omega)$ with scalar product
$$
	\langle u,v\rangle_{\H} = \int_U uv\,d\mu + \int_U\nabla u\cdot\nabla v \,d\mu.
$$
Moreover, let $\lambda>\f 1 2$ and define the sesquilinear forms
%\footnote{
%	$a(u,v)$ arises by applying Green's identity to $\int_\Om (-\Delta u) v \omega\,dx$.
%}
\begin{align}
	a^\alpha(u,v) &:= \int_U(\overline{\nabla u}\cdot\nabla v + \lambda \overline{u} v)\,d\mu + \int_U v\overline{\nabla u}\cdot\f{\nabla \omega}{\omega}\,d\mu + \alpha\int_{\del U}  \overline{u} v\,\omega\,dS &&\text{ on }\H,  \label{adef}  \\[0.5em]
	a^{\rm N}(u,v) &:= \int_U(\overline{\nabla u}\cdot\nabla v + \lambda \overline{u} v)\,d\mu + \int_U v\overline{\nabla u}\cdot\f{\nabla \omega}{\omega}\,d\mu  &&\text{ on }\H,\\[0.5em]
	a^{\rm D}(u,v) &:= \int_U(\overline{\nabla u}\cdot\nabla v + \lambda \overline{u} v)\,d\mu + \int_U v\overline{\nabla u}\cdot\f{\nabla \omega}{\omega}\,d\mu  &&\text{ on }\H_0.
\end{align}

\begin{lemma}\label{coercive}
	For $\lambda>\f 1 2$ and $\iota\in\{{\rm D}, {\rm N}, \alpha\},$ the form $a^\iota$ is continuous and coercive on $\H$ (on $\mathcal H_0$ in the case $\iota={\rm D}$).
\end{lemma}
\begin{proof}
	We will only treat the Robin case here, the other cases being analogous.
	Denote by $I$ the second term in \eqref{adef} and note that $\omega$ was chosen so that $|\nabla \omega| \leq \omega$. By H\"older's inequality with respect to $\mu$ one has
	\begin{align*}
		|I| \leq \underbrace{\left\| \dfrac{\nabla \omega}{ \omega} \right\|_\infty}_{\leq 1} \|\nabla u\|_{L^2(\mu)}\|v\|_{L^2(\mu)} 
		\leq \dfrac 1 2 \|\nabla u\|_{L^2(\mu)}^2 + \dfrac 1 2 \|v\|_{L^2(\mu)}^2,
	\end{align*}
	and thus
	\begin{align*}
		\bigl|a(u,u)\bigr| &\ge \|\nabla u\|_{L^2(\mu)}^2 + \lambda \|u\|_{L^2(\mu)}^2 + |\alpha|\bigl\|\omega^{\nicefrac 1 2}u\bigr\|_{L^2(\del U)}^2 + I \\[0.5em]
		&\geq \|\nabla u\|_{L^2(\mu)}^2 + \lambda \|u\|_{L^2(\mu)}^2 - \dfrac 1 2 \|\nabla u\|_{L^2(\mu)}^2 - \f 1 2 \|u\|_{L^2(\mu)}^2 \\[0.5em]
		&= \dfrac 1 2 \|\nabla u\|_{L^2(\mu)}^2 +\biggl(\lambda - \dfrac 1 2\biggr)\|u\|_{L^2(\mu)}^2,
	\end{align*}
	which shows coercivity in $\H$. Continuity follows by estimating the boundary term. By the trace theorem \cite[Prop. IX.18.1]{DiB} we have, for each $\delta>0,$
	\begin{align}
		\int_{\del U} |u|^2\omega\,dx &\leq 2\delta\|\nabla(\omega^{\nicefrac 1 2}u)\|_{L^2(U)}^2 + \f{C}{\delta}\|\omega^{\nicefrac{1}{2}}u\|_{L^2(U)}^2. 
		\label{onestar}
	\end{align}
	The first term can be estimated using the special choice of $\omega:$
	\begin{align}
		\|\nabla(\omega^{\nicefrac 1 2}u)\|_{L^2(U)}^2 &= \int_U \left|\omega^{\nicefrac 1 2} \nabla u + \f 1 2 u\f{\nabla\omega}{\omega^{\nicefrac 1 2}}\right|^2\,dx\nonumber\\[0.5em]
		&\leq 2\int_U \omega|\nabla u|^2\,dx + \frac{1}{2}\int_U |u|^2\f{|\nabla\omega|^2}{\omega}\,dx
		\nonumber\\[0.5em]
		&\leq 2\|\nabla u\|_{L^2(\mu)} + 2\left\| \f{\nabla\omega}{\omega} \right\|_{\infty}^2\int_U |u|^2\omega\,dx\nonumber \\[0.5em]
		&\leq 2\|\nabla u\|_{H^1(\mu)}^2.\label{twostar}
	\end{align}
	The desired continuity now follows immediately by combining (\ref{onestar}) and (\ref{twostar}). 
\end{proof}
\begin{lemma}\label{weightedsol}
	Let $f\in L^2(U)$, $\iota\in\{{\rm D}, {\rm N},\alpha\},$ and suppose that $\supp(f)$ is compact. Then the problem
	\begin{equation}\label{weightedeq}
		a^\iota(u,v) = \int_U \overline{f}v\,d\mu\quad\quad\forall v\in\H 
	\end{equation}
	has a solution in $\H$.
\end{lemma}
\begin{proof}
	By H\"older inequality, one has 
	$$\left|\int_U \overline{f}v\,d\mu\right| \leq \|f\|_{L^2(\mu)}\|v\|_{L^2(\mu)} \leq \|\omega\|_{L^\infty(\supp f)}\|f\|_{L^2(U)}\|v\|_{L^2(\mu)},$$
	 so $f\in \H'$. The assertion now follows from Lemma \ref{coercive} and the Lax-Milgram theorem for complex, non-symmetric sesquilinear forms \cite[Thm. VI.1.4]{TaylorLay}.
\end{proof}
\begin{proof}[Proof of Proposition \ref{main}] Again we focus on the Robin case, the other cases being analogous. Denote by $u$ the solution obtained from Lemma \ref{weightedsol}. Then $u\in H^1(U)$, since $\H\subset H^1(U)$. Moreover, let $\phi\in C_0^\infty(\R^d)$ be arbitrary and decompose it as $\phi=\omega\psi$. Then $\psi\in C_0^\infty(\R^d)\subset \H$ and one has
\begin{align*}
	\int_U\overline{\nabla u}\cdot\nabla\phi\,dx + \lambda\int_U \overline{u}\phi\,dx&+ \alpha\int_{\del U}\overline{u}\phi\,dS\\
	&= \int_U \overline{\nabla u}\cdot\big(  \omega\nabla\psi + \psi\nabla \omega \big)\,dx + \lambda\int_U \overline{u}\psi \omega\,dx + \alpha\int_{\del U}\overline{u}\psi\omega\,dS \\
	&= a^\alpha(u,\psi) \\
	&=\int_U \overline{f}\psi\,d\mu \\
	&=\int_U \overline{f}\phi\,dx.
\end{align*}
Thus, the function $u$ solves the problem
\begin{equation}
	\int_U\overline{\nabla u}\cdot\nabla\phi\,dx + \lambda\int_U \overline{u}\phi\,dx + \alpha\int_{\del U}\overline{u}\phi\,dS = \int_U \overline{f}\phi\,dx \quad\quad\forall \phi\in C_0^\infty(\R^d).
\end{equation}
Uniqueness of solutions and density of $C_0^\infty(\R^d)$ in $H^1(U)$ implies that $u$ is {the} weak solution in $H^1(U)$ to the Robin problem \eqref{GenRob}.

The estimates \eqref{est1}, \eqref{est2} follow from the coercivity of $a^\iota$.
\end{proof}

\section{Decomposition of the right-hand side}\label{decomposition}

In this section we consider the case of unbounded $\Om$. We conclude the proof of Theorem \ref{mainth} by decomposing the domain into cubes $Q_i,$ writing $f = \sum_i f\chi_{Q_i}$ and then applying the above results to each term $f\chi_{Q_i}.$ The following lemma shows uniform convergence with respect to the position of the cubes.
\begin{lemma}\label{boundedproof}
	Let $\eps_n\searrow0$ and $f_n\in L^2(\Om),$ $n\in{\mathbb N},$ be such that $\|f_n\|_{L^2(\Om)}\leq 1$ and $\supp (f_n)\subset Q_{i_n}$, where $Q_{i_n}=[0,1]^d+i_n$ with $ i_n\in\mathbb Z^d$. Let $u_n^\iota,\tilde u_n^\iota$ be the solutions to the problems 
	%($\iota\in\{{\rm D}, {\rm N}, \alpha\}$)
	\begin{align}\label{unshifted}
		A_{\eps_n}^\iota u_n^\iota = f_n|_{\Om_{\eps_n}}, \qquad
		A^\iota \tilde u_n^\iota = f_n,\quad\quad n\in{\mathbb N},\quad\quad\quad \iota\in\{{\rm D}, {\rm N}, \alpha\}.
	\end{align}
	Then 
	$\|J_{\eps_n} u_n^\iota-\tilde u_n^\iota\|_{L^2(\Om)}\to 0$ for all	 $\iota\in\{{\rm D}, {\rm N}, \alpha\}.$
\end{lemma}
\begin{proof}
 The idea of the proof is to use translation invariance, in order to shift $\supp(f_n)$ back near zero for every $n,$ and then use the Fr\'echet-Kolmogorov compactness criterion to obtain a convergent subsequence of $(J_{\eps_n}u_n^\iota-\tilde u_n^\iota)$; Theorem \ref{uniright} will identify its limit as zero. Since the following analysis is independent of the choice of boundary conditions, we henceforth omit $\iota$ to simplify notation.
	
	We now carry out the outlined strategy. We set, for \frank{$n\in{\mathbb N},$}
	\begin{align*}
		u_n^*(x) := u_n(x+i_n),\quad\quad
		%\\[0,2em]
		{\tilde u}_n^*(x) := \tilde u_n(x+i_n),\quad\quad
		%\\[0.2em]
		 f_n^*(x) := f_n(x+i_n).
	\end{align*}
	These functions still solve the problems (\ref{unshifted}) with $f_n$ replaced by ${f}_n^*$ and $\Omega$ replaced by $\Om-i_n$.
	The new sequence $f_n^*$ has the nice property that $\supp( f_n^*)\subset [0,1]^d$ for all $n$. In the following we consider $J_{\eps_n}u_n^*,\,\tilde u_n^*,\, f_n^*$ as elements of $L^2(\R^d)$ that are zero outside $\Om-i_n$.
	We will now show that ${\tilde u}_n^*- J_{\eps_n}u_n^*$ converges to zero in $L^2(\R^d)$. To this end, consider the bounded set
	\begin{equation}\label{Fprecompact}
		\mathcal F:=\{{\tilde u}_n^*-J_{\eps_n} u_n^* : n\in\mathbb N\}\subset L^2(\R^d).
	\end{equation}

	\begin{itemize}
		\item[\textit{Claim:}] 
		$\mathcal F$ is precompact in $L^2(\R^d)$.
	\end{itemize}	
	We postpone the proof of this claim to Lemma \ref{precompact}. We immediately obtain that $({\tilde u}_n^*-J_{\eps_n} u_n^*)$ has a convergent subsequence in $L^2(\R^d)$. Furthermore, Theorem \ref{uniright} shows that $\|{\tilde u}_n^*-J_{\eps_n} u_n^*\|_{L^2(K)}\to 0$ for every bounded $K\subset\R^d$ which identifies the limit of the subsequence as zero. 
	
	%Applying the above reasoning to 
	Arguing as above for all subsequences of $({\tilde u}_n^*-J_{\eps_n} u_n^*)$, we conclude that 
	%the convergence of
	 %the whole sequence
	${\tilde u}_n^*-J_{\eps_n} u_n^* \to 0$
		 %\quad\text{in }
		% holds 
		 in $L^2(\R^d).$ Since the shift $u\mapsto u(\cdot+i_n)$ is an isometry in $L^2(\R^d)$, this implies that $
		{\tilde u}_n-J_{\eps_n} u_n \to 0$ in $L^2(\Om).$
\end{proof}

\begin{lemma}\label{precompact}
	The set $\mathcal F$ defined in \eqref{Fprecompact} is precompact in $L^2(\R^d)$.
\end{lemma}

\begin{proof}
	We will use the notation and conventions from the previous proof and distinguish between the Dirichlet case and the Robin/Neumann cases. 
	%We start with the Dirichlet case.
	\paragraph{\it Dirichlet case.}
	Step 1: We have 
	$$
		\sup_n\bigl\|\tau_h({\tilde u}_n^*- J_{\eps_n}u_n^*)-({\tilde u}_n^*- J_{\eps_n}u_n^*)\bigr\|_{L^2(\R^d)}\to 0\quad{\rm as}\ h\to 0\quad\quad\forall n\in{\mathbb N},
	$$
	where $\tau_h$ denotes the operator of translation by $h.$ Indeed, the standard regularity theory implies 
	\begin{align*}
		\bigl\|\tau_h({\tilde u}_n^*- J_{\eps_n}u_n^*)-({\tilde u}_n^*-J_{\eps_n}u_n^*)\bigr\|_{L^2(\R^d)} &\leq\bigl\|\nabla ({\tilde u}_n^*-J_{\eps_n} u_n^*)\bigr\|_{L^2(\R^d)}|h| \\[0.4em]
			&\leq C\|f_n\|_{L^2(\Omega)}|h|.
	\end{align*}
	Step 2: Notice that
	\begin{align*}
		\sup_n\|{\tilde u}_n^*-J_{\eps_n} u_n^*\|_{L^2(\R^d\setminus B_R(0))} \to  0\quad{\rm as}\ R\to\infty,
	\end{align*}
	due to the following estimate in which we set $\omega_0(x):=\cosh(|x|)$.
	\begin{align*}
		\|{\tilde u}_n^*-J_{\eps_n} u_n^*\|^2_{L^2(\R^d\setminus B_R(0))} 
		\leq\; & 2\|\tilde u_n^*\omega_0\omega_0^{-1}\|_{L^2(\Om\setminus B_R(0))}^2 + 2\|J_\eps u_n^*\omega_0\omega_0^{-1}\|^2_{L^2((\R^d\setminus B_R(0))}\\[0.3em]
		\leq\; & 4M\|f_n^*\omega_0\|_{L^2(\R^d)}^2\|\omega_0^{-1}\|^2_{L^\infty(\R^d\setminus B_R(0))} \\[0.3em]	
		\stackrel{\mathmakebox[\widthof{=}]{\text{Prop. }\ref{main}}}{\qquad\quad\leq}\;&C\|f_n\|^2_{L^2(\Om)}\exp(-R).
	\end{align*}
	which completes Step 2. Applying the Fr\'echet-Kolmogorov theorem yields the precompactness of $\mathcal F$.
	
	\paragraph{\it Neumann and Robin case.}
	
	Here the strategy is the same, but matters are complicated by the fact that $J_{\eps_n}u_n^*$ is not in $H^1(\R^d)$. To show that $\mathcal F$ is precompact, 
	%we will show that $\mathcal F = \mathcal F_1 + \mathcal F_2$, where $\mathcal F_1,\mathcal F_2$ are precompact in $L^2(\R^d)$. More precisely, 
	we decompose elements in ${\mathcal F}$ as 
	\begin{align*}
		{\tilde u}_n^*- J_{\eps_n}u_n^* = ({\tilde u}_n^*-\mathcal T_{\eps_n} u_n^*) + (\mathcal T_{\eps_n}- J_{\eps_n})u_n^*,
	\end{align*}
define $\mathcal F_1:=\{{\tilde u}_n^*-\mathcal T_{\eps_n} u_n^*:n\in\mathbb N\}$, $\mathcal F_2:=\{(\mathcal T_{\eps_n} - J_{\eps_n})u_n^*:n\in\mathbb N\}$ and show that $\mathcal F_1$ and $\mathcal F_2$ are precompact in $L^2(\R^d).$ We will begin by showing that $\mathcal F_1$ is precompact. 
\frank{
To this end, denote by $\mathcal E:H^1(\Om)\to H^1(\R^d)$ an extension operator satisfying $\mathcal Eu|_\Om=u$ and $\|\mathcal Eu\|_{H^1(\R^d)}\leq C\|u\|_{H^1(\Om)}$ for all $u\in H^1(\Om)$ (cf. \cite[Theorem 5.24]{Adams}). 

Clearly, for every $\xi\in\R^d$ the operators $\mathcal E_\xi:H^1(\Om-\xi)\to H^1(\R^d)$ defined by $\mathcal E_\xi u:=\tau_{\xi}\mathcal E\tau_{-\xi} u$ satisfy $\|\mathcal E_\xi\|_{\mathcal L(H^1(\Om-\xi),H^1(\R^d))} = \|\mathcal E\|_{\mathcal L(H^1(\Om),H^1(\R^d))} $.
We start by proving that 
	\begin{align*}
		\sup_n\big\|\tau_h\mathcal E_{i_n}(\tilde u_n^*-\mathcal T_{\eps_n} u_n^*) - \mathcal E_{i_n}(\tilde u_n^*-\mathcal T_{\eps_n} u_n^*)\big\|_2\to 0\quad{\rm as}\ h\to0
	\end{align*}
	This readily follows from the estimate
	\begin{align*}
		\bigl\|\tau_h\mathcal E_{i_n}({\tilde u}_n^*-\mathcal T_{\eps_n} u_n^*)-\mathcal E_{i_n}({\tilde u}_n^*-\mathcal T_{\eps_n} u_n^*)\bigr\|_{L^2(\R^d)}
		&\leq \bigl\|\nabla \mathcal E_{i_n}({\tilde u}_n^*-\mathcal T_{\eps_n} u_n^*)\bigr\|_{L^2(\R^d)}|h| \\[0.4em]
		&\leq C\|{\tilde u}_n^*-\mathcal T_{\eps_n} u_n^*\|_{H^1(\Om+i_n)}|h|\\[0.4em]
		&\leq C\| J_{\eps_n}f_n^*\|_{L^2(\Omega+i_n)}|h|\\[0.4em]
		&\leq C|h|.
	\end{align*}
	Next we prove that
	\begin{align*}
		\sup_n\bigl\|\mathcal E_{i_n}({\tilde u}_n^* - \mathcal T_{\eps_n} u_n^*)\bigr\|_{L^2(\R^d\setminus B_R(0))} \to  0\quad{\rm as}\ R\to\infty.
	\end{align*}
 Indeed, notice first that
 % have the estimate
	\begin{align}
		\bigl\|\mathcal E_{i_n}({\tilde u}_n^*-\mathcal T_{\eps_n} u_n^*)\bigr\|^2_{L^2(\R^d\setminus B_R(0))}   
		&\leq	  C\left( \|{\tilde u}_n^*\|^2_{L^2((\Om+i_n)\setminus B_R(0))} + \| \mathcal T_{\eps_n}u_n^*\|^2_{L^2((\Om_{\eps_n}+i_n)\setminus B_R(0))} \right)\nonumber\\[0.4em]	
		&=   C\left( \|{\tilde u}_n\|^2_{L^2(\Om\setminus B_R(i_n))} + \| \mathcal T_{\eps_n}u_n\|^2_{L^2((\Om_{\eps_n})\setminus B_R(i_n))}\right),\label{T_est}
	\end{align}
	} 
	To treat the two terms on the right-hand side we apply Lemma \ref{harmext} (ii) and Proposition \ref{main} with $\omega_{i_n}(x) = \cosh(|x-i_n|)$ as follows. For the second term in (\ref{T_est}), we obtain 
	\begin{align*}
		\| \mathcal T_{\eps_n}u_n\|_{L^2(\Om_{\eps_n}\setminus B_R(i_n))} 
		&\leq C\big(\|u_n\|_{L^2(\Om\setminus B_{\nicefrac R 2}(i_n))}+\|\nabla u_n\|_{L^2(\Om\setminus B_{\nicefrac R 2}(i_n))}\big) \\[0.4em]
		&\leq \bigl\|\omega_{i_n}^{\nicefrac 1 2}\omega_{i_n}^{-\nicefrac 1 2} u_n\bigr\|_{L^2(\Om\setminus B_{\nicefrac R 2}(i_n))} + \bigl\|\omega_{i_n}^{\nicefrac 1 2}\omega_{i_n}^{-\nicefrac 1 2} \nabla u_n\bigr\|_{L^2(\Om\setminus B_{\nicefrac R 2}(i_n))} \\[0.4em]
		&\leq  C\Big(\bigl\|\omega_{i_n}^{\nicefrac 1 2} u_n\bigr\|_{L^2(\Om\setminus B_{\nicefrac R 2}(i_n))} \\
		&\hspace{3cm} + \bigl\|\omega_n^{\nicefrac 1 2}\nabla u_n\bigr\|_{L^2(\Om\setminus B_{\nicefrac R 2}(i_n))}\Big)\|\omega_{i_n}^{-\nicefrac 1 2}\|_{L^\infty(\Om\setminus B_{\nicefrac R 2}(i_n))} \\[0.4em]
		&\leq CM\bigl\|f_n\omega_{i_n}^{\nicefrac 12}\bigr\|_{L^2(\Om)}\exp(-R/3) \\[0.4em]
		&\leq 2CM\exp(-R/3),
	\end{align*}
	
	where we used the fact that $\omega_{i_n}$ is bounded by 2 on $\supp f_n$. With an analogous calculation for the first term in (\ref{T_est}),
	%$\tilde u_n$ 
	we finally find
	\begin{equation}
		\bigl\|\mathcal E_{i_n}({\tilde u}_n^*-\mathcal T_{\eps_n} u_n^*)\bigr\|_{L^2(\R^d\setminus B_R(0))} \leq C\exp(-R/3),
	\end{equation}
	with $C$ independent of $n$. Applying the Fr\'echet-Kolmogorov theorem yields the precompactness of the set $\{\mathcal E_{i_n}({\tilde u}_n^*-\mathcal T_{\eps_n} u_n^*):n\in\mathbb N\}$. Finally, noting that $\mathcal F_1 = \{\mathcal E_{i_n}({\tilde u}_n^*-\mathcal T_{\eps_n} u_n^*):n\in\mathbb N\}\!\cdot\!\chi_{\Om}$ and that multiplication by $\chi_\Om$ is a bounded operator on $L^2(\R^d)$ we obtain precompactness of $\mathcal F_1$.
	\\\\
	To prove precompactness of $\mathcal F_2$, first note that by Lemma \ref{harmext} (iii) for any $\delta>0$ there exists a $n_0$ such that 
	\begin{align*}
		\bigl\|(J_{\eps_n} - \mathcal T_{\eps_n})u_n^*\bigr\|_2 < \delta \quad\forall n>n_0.
	\end{align*}
	Let us fix  arbitrary $\delta>0$ and $n_0$ as above. It  remains to estimate the terms 
	\[
	\big\| \tau_h(J_{\eps_n} - \mathcal T_{\eps_n})u_n^* - (J_{\eps_n} - \mathcal T_{\eps_n})u_n^*\big\|_{L^2({\mathbb R}^d)},\quad\quad n\leq n_0,
	\]
	 but these are only finitely many, which clearly converge to zero individually as $h\to 0$, and hence 
	$$\sup_{n\leq n_0}\big\| \tau_h(J_{\eps_n} - \mathcal T_{\eps_n})u_n^* - (J_{\eps_n} - \mathcal T_{\eps_n})u_n^*\big\|_2\to 0\quad\text{ as }h\to 0$$
	Altogether we have shown that 
	\begin{align*}
		\sup_n \big\| \tau_h(J_{\eps_n} - \mathcal T_{\eps_n})u_n^* - (J_{\eps_n} &- \mathcal T_{\eps_n})u_n^*\big\|_{L^2({\mathbb R}^d)}\\
		&\leq \max\left\{ \sup_{n\leq n_0} \big\| \tau_h(J_{\eps_n} - \mathcal T_{\eps_n})u_n^* - (J_{\eps_n} - \mathcal T_{\eps_n})u_n^*\big\|_2\, ,\, 2\delta \right\} \\
		&\!\xrightarrow{h\to 0}\; 2\delta.
	\end{align*}
	Since $\delta>0$ was arbitrary we finally get
	\begin{align*}
		\lim_{h\to 0}\,\sup_{n\in\mathbb N} \big\| \tau_h(J_{\eps_n} - \mathcal T_{\eps_n})u_n^* - (J_{\eps_n} - \mathcal T_{\eps_n})u_n^*\big\|_{L^2({\mathbb R}^d)} = 0.
	\end{align*}
	This completes the first Fr\'echet-Kolmogorov-condition. 
	The proof of the second condition
	\begin{align*}
		\sup_n\bigl\|(J_{\eps_n}-\mathcal T_{\eps_n}) u_n^*\bigr\|_{L^2(\R^d\setminus B_R(0))} \to  0\quad{\rm as}\ R\to\infty
	\end{align*}
	is analogous to the case of $\mathcal F_1$.
	Applying the Fr\'echet-Kolmogorov theorem yields  precompactness of $\mathcal F_2$ and completes the proof.
\end{proof}

\begin{corollary}\label{Cor}
	There exists $\delta_\eps$ with $\delta_\eps\stackrel{\scriptscriptstyle{\eps\to 0}}{\longrightarrow} 0$ such that 
	$$
		\bigl\|(J_\eps (A^\iota)^{-1}-(A_\eps^\iota)^{-1}J_\eps )(f\chi_{Q_i\cap\Om_\eps})\bigr\|_{L^2(\Omega)} \leq \delta_\eps\|f\chi_{Q_i}\|_{L^2(\Omega)}
	$$
	for all $f\in L^2(\Om)$ and $i\in\mathbb Z^d$.
\end{corollary}
\begin{proof}
	We argue by contradiction. Suppose that there is no such function $\delta_\eps$. Then there exist sequences $\eps_n,f_n,i_n$ with $\|f_n\|_{L^2(\Omega)}= 1$ such that $\|(J_\eps (A^\iota)^{-1}-(A_{\eps_n}^\iota)^{-1}J_\eps )(f_n\chi_{Q_{i_n}\cap\Om_{\eps_n}})\|_{L^2(\Omega)}$ does not converge to zero, which is a contradiction to Lemma \ref{boundedproof}.
\end{proof}
In order to finalise the decomposition, we require he following two lemmas.
\begin{lemma}\label{lemma2}
	Suppose that $f\in L^2(\Omega),$ and denote 
	\[
	u_i:=\bigl(J_\eps (A^\iota)^{-1}-(A_\eps^\iota)^{-1}J_\eps\bigr)(f\chi_{Q_i\cap\Om_\eps}),\quad i\in{\mathbb Z}^d.
	\]
	Then one has 
	\begin{equation}
		\bigl|\langle u_i,u_j\rangle_{L^2(\Omega)}\bigr|\leq C\|f\chi_{Q_i}\|_{L^2(\Omega)}\|f\chi_{Q_j}\|_{L^2(\Omega)}\,\exp(-|i-j|/2)
	\end{equation}
	for all $i,j\in{\mathbb Z}^d$ with $i\neq j$, where $\langle\cdot,\cdot\rangle_{L^2(\Omega)}$ denotes the standard inner product in $L^2(\Omega).$
\end{lemma}
\begin{proof}
	For convenience we write $f_i:=f\chi_{Q_i},$ $i\in{\mathbb Z}^d$. Denote $\omega_i(x)=\cosh(|x-i|)$ and note that by Proposition \ref{main} we have $\|\omega_i^{\nicefrac 1 2}u_i\|_{L^2(\Om)}\leq C \|f_i\omega_i^{\nicefrac 1 2}\|_{L^2(\Om)}$.
	The statement of the lemma is a consequence of the following estimate:
	\begin{align*}
		\bigl|\langle u_i,u_j\rangle_{L^2(\Omega)}\bigr| &\leq \int_{\Om} |u_i(x)||u_j(x)|\,dx\nonumber\\[0.4em]
		&= \int_{\Om} \bigl(|u_i(x)|\omega_i^{\nicefrac 1 2}\bigr)\big(|u_j(x)|\omega_j^{\nicefrac 1 2}\big)\omega_i^{-\nicefrac 1 2}\omega_j^{-\nicefrac 1 2}\,dx\nonumber\\[0.4em]
		&\leq \bigl\|u_i\omega_i^{\nicefrac 1 2}\bigr\|_{L^2(\Om)}\bigl\|u_j\omega_j^{\nicefrac 1 2}\bigr\|_{L^2(\Om)}\bigl\|\omega_i^{-\nicefrac 1 2}\omega_j^{-\nicefrac 1 2}\bigr\|_{L^\infty(\Om)}\\[0.6em]
		&\leq C\|f_i\omega_i^{\nicefrac 1 2}\|_{L^2(\Om)}\,\|f_j\omega_j^{\nicefrac 1 2}\|_{L^2(\Om)}\,\|\omega_0^{-\nicefrac 1 2}\omega_{j-i}^{-\nicefrac 1 2}\|_{L^\infty(\Om)}\\[0.6em]
		&\leq C\|f_i\|_{L^2(\Om)}\|f_j\|_{L^2(\Om)}\exp\bigl(-|i-j|/2\bigr),
	\end{align*}
	where we use the fact that $\supp( f_i)\subset Q_i$ and $\omega_i|_{Q_i}\leq 2$. 
\end{proof}

\begin{lemma}	\label{inequality}
	Suppose that $f\in C_0^\infty(\Om_\eps)$ and define $u_i:=(J_\eps (A_\eps^\iota)^{-1}-(A^\iota)^{-1}J_\eps )(f\chi_{Q_i}),$ $i\in{\mathbb Z}^d.$ Then for every $n>1$ one has the inequality
	\begin{equation}\label{ConIneq}
		\left\|  \sum_{m=1}^N u_{i_m}  \right\|_{L^2(\Om)}^2 \leq  C\left( n^3\sum_{m=1}^N \|u_{i_m}\|_{L^2(\Om)}^2 + \|f\|_{L^2(\Om_\eps)}
		\exp(-n/3)\right),
	\end{equation}
	where $N$ is the number of cubes such that $Q_{i_k}\cap\supp( f)\neq\emptyset,$ and $C,n$ do not depend on $N$.
\end{lemma}
\begin{proof}
	\begin{align}
		\left\|  \sum_{m=1}^N u_{i_m}  \right\|_{L^2(\Omega)}^2 
		&\leq \sum_{m,p=1}^N\langle u_{i_m},u_{j_p}\rangle_{L^2(\Omega)}\nonumber \\
		&= \sum_{k=0}^\infty \left(\sum_{|i-j|\in[k,k+1)}\langle u_i,u_j\rangle_{L^2(\Omega)}\right)\nonumber\\
		&\leq \sum_{k=0}^n \left(\sum_{|i-j|\in[k,k+1)}\|u_i\|_{L^2(\Omega)}\|u_j\|_{L^2(\Omega)}\right) + \sum_{k=n}^\infty \left(\sum_{|i-j|\in[k,k+1)}\langle u_i,u_j\rangle_{L^2(\Omega)}\right)\nonumber\\
		&\leq \sum_{k=0}^n\; \sum_{|i-j|\in[k,k+1)}\left( \f{\|u_i\|_{L^2(\Omega)}^2}{2}+\f{\|u_j\|_{L^2(\Omega)}^2}{2} \right) + \sum_{k=n}^\infty \left(\sum_{|i-j|\in[k,k+1)}\langle u_i,u_j\rangle_{L^2(\Omega)}\right)\nonumber\\
		&\leq \sum_{k=0}^n \;\sum_{m=1}^N\left(\sum_{\{j:|i_m-j|\in[k,k+1)\}}\|u_{i_m}\|_{L^2(\Omega)}^2\right) + \sum_{k=n}^\infty \left(\sum_{|i-j|\in[k,k+1)}\langle u_i,u_j\rangle_{L^2(\Omega)}\right)\nonumber\\[0.4em]
		&\leq C \sum_{k=1}^n k^2\sum_{m=1}^N\|u_{i_m}\|_{L^2(\Omega)}^2 + \sum_{k=n}^\infty \left(\sum_{|i-j|\in[k,k+1)}\langle u_i,u_j\rangle_{L^2(\Omega)}\right)\nonumber\\[0.3em]
		&\leq C n^3\sum_{m=1}^N\|u_{i_m}\|_{L^2(\Omega)}^2 + \sum_{k=n}^\infty \left(\sum_{|i-j|\in[k,k+1)}\langle u_i,u_j\rangle_{L^2(\Omega)}\right).\label{last_expr}
	\end{align}
	We now study the last term of (\ref{last_expr}) . It follows from Lemma \ref{lemma2} that 
	\[
	\bigl|\langle u_i,u_j\rangle_{L^2(\Omega)}\bigr|\leq C\|f_i\|_{L^2(\Omega)}\|f_j\|_{L^2(\Omega)}e^{-\f 1 2|i-j|}.
	\] 
	Using this fact and fixing $k$ for the moment, we obtain
	\begin{align*}
		\left| \sum_{|i-j|\in[k,k+1)}\langle u_i,u_j\rangle_{L^2(\Omega)}\right|  
			&\leq  C\sum_{|i-j|\in[k,k+1)} \|f_i\|_{L^2(\Omega)}\|f_j\|_{L^2(\Omega)}\exp(-|i-j|/2) \\
			&\leq C\sum_{|i-j|\in[k,k+1)} \left( \f{\|f_i\|_{L^2(\Omega)}^2}{2}+\f{\|f_j\|_{L^2(\Omega)}^2}{2} \right)\exp(-|i-j|/2) \\
			&\leq C\sum_{m=1}^N \|f_{i_m}\|_{L^2(\Omega)}^2 k^2\exp(-k/2) \\[0.4em]
			&= C\|f\|_{L^2(\Omega)}^2k^2\exp(-k/2) \\[0.6em]
			&\leq C\|f\|_{L^2(\Omega)}^2\exp(-k/2) .
	\end{align*}
	Summing this inequality from $k=n$ to infinity concludes the proof.
\end{proof}

Combining the above lemmas, we have the following quantitative statement.
\begin{prop}
	Suppose that $f\in C_0^\infty(\Om_\eps)$. Then for every $n\in\mathbb N$,
	$$\bigl\|(J_\eps (A_\eps^\iota)^{-1}-(A^\iota)^{-1}J_\eps)f\bigr\|_{L^2(\Omega)}^2\leq C\big(n^3 \delta_\eps^2 +\exp(-n/3)\big)\|f\|_{L^2(\Omega)}^2$$
	for some $C>0$, where $\delta_\eps$ was defined in Corollary \ref{Cor}.
\end{prop}
\begin{proof}
We denote $u_i^\eps:=(J_\eps (A_\eps^\iota)^{-1}-(A^\iota)^{-1}J_\eps)(f\chi_{Q_i}),$ $i\in{\mathbb R}^d,$ and estimate
\begin{align*}
	\bigl\|(J_\eps (A_\eps^\iota)^{-1}-(A^\iota)^{-1}J_\eps)f\bigr\|_{L^2(\Omega)}^2 
	=\;&  \left\|\sum_{m=1}^N u_{i_m}^\eps\right\|_{L^2(\Omega)}^2 \\
	\stackrel{\mathmakebox[\widthof{=}]{\text{Lemma \ref{inequality}}}}{\hspace{2cm}\leq}\; & C\left( n^3\sum_{m=1}^N \|u_{i_m}^\eps\|_{L^2(\Omega)}^2 + e^{-\nicefrac{n}{3}}\|f\|_{L^2(\Omega_\eps)}\right) \\
	\stackrel{\mathmakebox[\widthof{=}]{\text{Cor. \ref{Cor}}}}{\hspace{1.6cm}\leq}\; & C\left( n^3\delta_\eps^2\sum_{m=1}^N \|f_{i_m}\|_{L^2(\Omega_\eps)}^2 + e^{-\nicefrac{n}{3}}\|f\|_{L^2(\Omega_\eps)}\right) \\[0.5em]
	=\;& C\left( n^3\delta_\eps^2 + e^{-\nicefrac{n}{3}}\right)\|f\|_{L^2(\Omega)}^2.
	\end{align*}
\end{proof}

\paragraph{Proof of Theorem \ref{mainth}.}
	Let $g\in L^2(\Om_\eps)$ with $\|g\|_{L^2(\Omega_\eps)}\leq 1$. Fix $\delta>0$ and choose $f\in C_0^\infty(\Om_\eps)$ such that $\|g-f\|_{L^2(\Omega_\eps)}^2<\delta$ and choose $n\in\mathbb N$ such that $\exp(-n/3)\leq\delta$. Now compute
	\begin{align*}
		\bigl\|(J_\eps (A_\eps^\iota)^{-1}-(A^\iota)^{-1}J_\eps)g\bigr\|_{L^2(\Omega)}^2 &\leq 2 \bigl\|(J_\eps (A_\eps^\iota)^{-1}-(A^\iota)^{-1}J_\eps)f\bigr\|_{L^2(\Omega)}^2 + 2\bigl\|(J_\eps (A_\eps^\iota)^{-1}-(A^\iota)^{-1}J_\eps)(g-f)\bigr\|_{L^2(\Omega)}^2\\[0.5em]
		&\leq C\Big(\big(n^3\delta_\eps^2 + \exp(-n/3)\big)\|f\|_{L^2(\Omega_\eps)}^2 + \underbrace{\bigl\|J_\eps (A_\eps^\iota)^{-1}\!-\!(A^\iota)^{-1}J_\eps\bigr\|^2}_{\text{bounded}}\|g-f\|_{L^2(\Omega_\eps)}^2\Big)\\
		&\leq C(n^3\delta_\eps^2+\delta)\|g\|_{L^2(\Omega_\eps)}^2 + C\delta,
		%\Rightarrow\quad\; 
	\end{align*}
		hence
		\[
		\sup_{\|g\|_{L^2(\Omega_\eps)}\leq 1}\bigl\|(J_\eps (A_\eps^\iota)^{-1}-(A^\iota)^{-1}J_\eps )g\bigr\|_{L^2(\Omega)}^2 
		\leq Cn^3\delta_\eps^2+C\delta + C\delta,
		\]
		and therefore
		\[
		\limsup_{\eps\to 0}\bigl\|(J_\eps (A_\eps^\iota)^{-1}-(A^\iota)^{-1}J_\eps)\bigr\|_{\mathcal L(L^2(\Om_\eps),L^2(\Om))}^2 
		\leq C\delta.
		\]
	Since $\delta>0$ is arbitrary, the result follows.\qed

\section{Behaviour of the Semigroup}\label{semigroup}

In this section we want to give an application of Theorem \ref{mainth}. In particular, we focus on the non-selfadjoint operator $A_\alpha$ and study the large-time behaviour of its semigroup. In order to do this, we shall first study the numerical range of the Robin Laplacians more closely. In the remainder of this section, unless otherwise stated, the symbols $\|\cdot\|$ and $\langle\cdot,\cdot\rangle$ will denote the $L^2$ (operator-) norm and scalar product, respectively, and the symbol $\Sigma_{\theta}$ denotes a sector of half-angle $\theta$ in the complex plane.

\subsection{Decay of $\boldsymbol{{\rm e}^{-t(A^\alpha-\text{Id})}}$}

Let $\alpha\in\mathbb C$ and assume $\re\alpha>0$. We want to study the decay properties of the heat semigroup $e^{t(\Delta-\mu_\alpha)}$. To this end, let us denote by $B^\alpha:=A^\alpha-\text{Id}$ the Robin Laplacian on $\Om$. It is our goal to derive estimates on the numerical range of $B^\alpha$. Let $u\in\mathcal D(B^\alpha) = \mathcal D(A^\alpha)$ and assume that $\|u\|_{L^2(\Om)}=1$. Notice that
\begin{align*}
	\langle B^\alpha u,u\rangle &= \int_\Om |\nabla u|^2\,dx + \mu_\alpha \int_{\Om} |u|^2\,dx + \alpha\int_{\del\Omega}|u|^2\,dS\\[0.5em]
	&= \|\nabla u\|^2 + \mu_\alpha + \alpha\|u\|_{L^2(\del\Om)}^2,
\end{align*}
%Therefore, we have
and therefore 
\begin{align*}
	\re\bigl\langle B^\alpha u,u\bigr\rangle &\geq \re\mu_\alpha + \re\alpha \|u\|_{L^2(\del\Om)}^2, \\[0.5em]
	|\im\langle B^\alpha u,u\rangle| &\leq |\im\mu_\alpha| + |\im\alpha|\|u\|_{L^2(\del\Om)}^2.
\end{align*}
Now, let $\lambda\in(0,\re\mu_\alpha)$ and compute
\begin{align}
	\bigl|\im\bigl\langle (B^\alpha-\lambda) u,u\bigr\rangle\bigr| &\leq |\im\mu_\alpha| + |\im\alpha|\|u\|_{L^2(\del\Om)}^2\nonumber\\[0.5em]
	&= \f{|\im\mu_\alpha|}{\re\mu_\alpha}\re\mu_\alpha + \f{|\im\alpha|}{\re\alpha}\re\alpha\|u\|_{L^2(\del\Om)}^2.\label{B_est}
\end{align}
Recall from Section \ref{intro} that $\mu_\alpha =\alpha S_d/2^d$ and hence $|\im\mu_\alpha|/\re\mu_\alpha=|\im\alpha|/\re\alpha.$ Combining this with (\ref{B_est}), we obtain
\begin{align*}
	\bigl|\im\langle (B^\alpha-\lambda) u,u\rangle\bigr| &\leq \f{|\im\alpha|}{\re\alpha}\left(\re\mu_\alpha + \re\alpha\|u\|_{L^2(\del\Om)}^2\right) \\
	&\leq \f{|\im\alpha|}{\re\alpha}\bigl(\re\bigl\langle(B^\alpha-\lambda)u, u\bigr\rangle + \lambda \bigr)\\[0.5em]
	&\leq \f{|\im\alpha|}{\re\alpha-\dfrac{\lambda}{2^{-d}S_d}}\re\bigl\langle(B^\alpha-\lambda)u,u\bigr\rangle.
\end{align*}
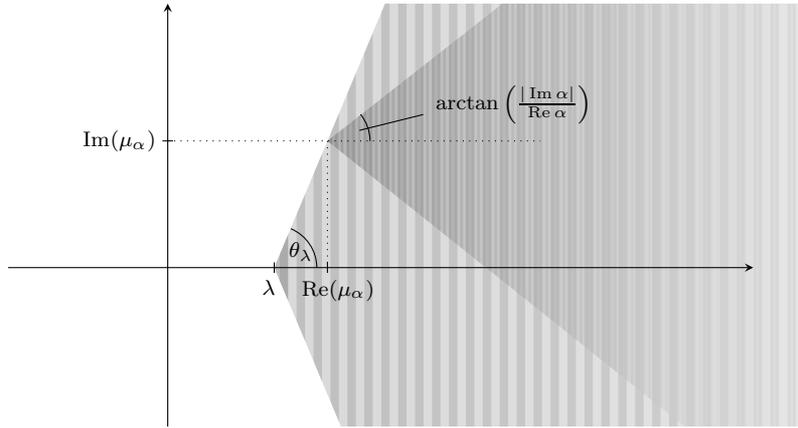
\begin{figure}[htbp]
	\centering
	%\documentclass[11pt,a4paper]{article} 
%
%\usepackage{etex}
%\usepackage[english]{babel}
%\usepackage{tikz}
%\usepackage{pgfplots}
%\usetikzlibrary{fadings,fpu,arrows,decorations.pathmorphing,backgrounds,positioning,fit,petri,patterns}
%\usepackage[utf8]{inputenc}
%\usepackage{graphicx}
%\usepackage[colorlinks=true,linkcolor=black,citecolor=black]{hyperref}
%\usepackage{amsfonts}
%\usepackage{enumerate}
%\usepackage{booktabs}
%\usepackage{array} 
%\usepackage{paralist} 
%\usepackage{subfig} 
%\usepackage{mathtools}
%\usepackage{tabu}
%\usepackage{amsthm}
%\usepackage{empheq}
%\usepackage{amsopn}
%\usepackage{dsfont}
%\usepackage{wrapfig}
%\usepackage[font=small]{caption}
%
%
%\begin{document}
\begin{tikzpicture}[>=stealth, scale=0.7]
    \clip (-6,-3) rectangle (9,5);
	\shade[left color=gray, opacity=0.3, right color=white]
    (-1,0) -- (7,-19.2) arc (-100:100:19.5) -- (-1,0);
	\shade[left color=gray, opacity=0.4, right color=white]
    (0,2.4) -- (9,-4.85) arc (-20:20:21) -- (0,2.4);
	
%	\draw (1,5) -- (1,0);
%	\draw (0.8,3.2) -- (1.2,3.2);
%	\draw (0.8,1.6) -- (1.2,1.6);
	
	\draw[dotted] (0,2.4) -- (0,0);
	\draw[dotted] (-3,2.4) -- (4,2.4);
	\draw (0.8,2.4) arc (0:39:0.8);
	
	\draw (0.6,2.6) -- (1.8,2.9);
	\draw (3.5,3.1) node{\footnotesize$\arctan\Big(\tfrac{|\operatorname{ Im}\alpha|}{\operatorname{ Re}\alpha}\Big)$};
	\draw (-0.2,0) arc (0:67:0.8);
	\draw (-0.5,0.3) node{\footnotesize$\theta_\lambda$};
	
	\draw[->] (-6,0) -- (8,0);
	\draw[->] (-3,-5) -- (-3,5);
	
	\draw (0,-0.1) -- (0,0.1);
	\draw (0.2,-0.05) node[below]{\footnotesize${\rm Re}(\mu_\alpha)$};
	\draw (-3.1,2.4) -- (-2.9,2.4);
	\draw (-3.05,2.4) node[left]{\footnotesize${\rm Im}(\mu_\alpha)$};
	
	\draw (-1,-0.1) -- (-1,0.1);
	\draw (-1.1,-0.05) node[below]{\footnotesize$\lambda$};
\end{tikzpicture}
%\end{document}
	\caption{The sector of decay and angle $\theta_\lambda$ for $B^\alpha$.}
\end{figure}
Using standard generation theorems about analytic semigroups, the next statement follows.
%, have shown the following
\begin{prop}
	The operator $-(B^\alpha-\lambda)$ generates a bounded analytic semigroup in the sector $\Sigma_{\f{\pi}{2}-\theta_\lambda}$, where 
	\[
	\theta_\lambda= \arctan\left(\f{|\im\alpha|}{\re\alpha-\dfrac{\lambda}{2^{-d}S_d}}\right).
	\]
	 Equivalently, $-B^\alpha$ generates an analytic semigroup with
	$$\bigl\|\exp(-zB^\alpha)\bigr\|\leq\exp(-\lambda z)\qquad\forall z\in \Sigma_{\f{\pi}{2}-\theta_\lambda}.$$
\end{prop}
\begin{proof}
	See \cite[Ch. IX.1.6]{Kato}.
\end{proof}

\subsection{Decay of $\boldsymbol{{\rm e}^{-t(A^\alpha_\eps-\text{Id})}}$}

In this section we denote $B^\alpha_\eps := A^\alpha_\eps-\text{Id}$. By calculations analogous to the above, we have
\begin{align*}
	\bigl|\im\langle B^\alpha_\eps  u,u\rangle\bigr| &\leq 
	 \f{|\im\alpha|}{\re\alpha}\re\langle B^\alpha_\eps  u,u \rangle,
\end{align*}
that is, $B^\alpha_\eps$ is sectorial with sector $\Sigma_{\theta_0},$ where 
$\theta_0 = \arctan(|\im\alpha|/\re\alpha),$ and hence generates a bounded analytic semigroup in the sector $\Sigma_{\f{\pi}{2}-\theta_0}$. In this subsection we improve this {\it a priori} result using spectral convergence. To this end, let $\delta>0$ and define the compact set
$$K_\delta:=\bigl\{ x+iy : x\in[0,\re\mu_\alpha],\,y\in\bigl[-|\im\mu_\alpha|,|\im\mu_\alpha|\bigr]\bigr\}.$$
Note that then $\Sigma_{\theta_0}\cap \{\re z\leq \re\mu_\alpha-\delta\}\subset K_\delta$. By \cite[Th. III.2.3]{EE} one has $K_\delta\subset\rho(B^\alpha)$ for every $\delta>0$. Applying Corollary \ref{SpCon} we see that for every $\delta>0$ there exists a $\eps_0>0$ such that $K_\delta\subset\rho(B^\alpha_\eps)$ for all $\eps<\eps_0$.

In particular we have shown that the resolvent norm $\|(B^\alpha_\eps-z)^{-1}\|$ is bounded on $\Sigma_{\theta_0}\cap \{\re z\leq \re\mu_\alpha-\delta\}$. By a trivial calculation analogous to the previous subsection this leads to the following statement.
\begin{lemma}\label{thetalambdadelta}
	For every $\lambda\in(0,\re\mu_\alpha-\delta)$ one has 
	$$\sigma(B^\alpha_\eps-\lambda)\subset \Sigma_{\theta_\lambda^\delta},\quad\quad \theta_\lambda^\delta= \arctan\biggl( \f{|\im\mu_\alpha|}{\re\mu_\alpha-\lambda-\delta} \biggr).$$
\end{lemma} 
Furthermore, we obtain the following lemma. 
%combine our results to obtain
\begin{lemma}
	For every $\lambda\in(0,\re\mu_\alpha-\delta)$ one has $\mathbb C\setminus\Sigma_{\theta_\lambda^\delta}\subset\rho(B^\alpha_\eps-\lambda)$ and there exists a $M=M(\lambda,\delta)>0$ such that 
	$$
		\bigl\|(B^\alpha_\eps-\lambda-z)^{-1}\bigr\| \leq \f{M}{|z|} \qquad\forall z\in \mathbb C\setminus\Sigma_{\theta_\lambda^\delta}.
	$$
\end{lemma}
\begin{proof}
	This is obtained by combining Lemma \ref{thetalambdadelta} with the following two facts:
	%\begin{itemize}
		%\item[(i)] 
		$$
		|\im\langle B^\alpha_\eps  u,u\rangle| \leq 
	 	\f{|\im\alpha|}{\re\alpha}\re\langle B^\alpha_\eps  u,u \rangle,\qquad\qquad
		%\item[(ii)]
	 	\bigl\|(B^\alpha_\eps-z)^{-1}\bigr\| \leq C\quad\text{ on } \;K_\delta.
		$$
	%\end{itemize}
\end{proof}

By the theory of analytic semigroups ({\it cf.} \cite[Ch. IX.1.6]{Kato}), we immediately obtain the following corollary.
\begin{corollary}
	For all $\lambda\in(0,\re\mu_\alpha-\delta),$ the operator $B^\alpha_\eps-\lambda $ generates a bounded analytic semigroup in the sector $\Sigma_{\f{\pi}{2}-\theta_\lambda^\delta}$.
\end{corollary}
This yields the main result of this section, as follows.
\begin{theorem}
	For every $\delta>0$ there exists $\eps_0>0$ such that for every $\lambda\in(0,\re\mu-\delta)$ there exists $M>0$ such that
	$$
		\bigl\|\exp(-zB^\alpha_\eps)\bigr\| \leq M\exp(-\lambda\re z)\quad\quad\forall z\in\Sigma_{\theta_\lambda^\delta},\ \ \eps\in(0,\eps_0).
	$$
	%for all $\eps\in(0,\eps_0)$.
\end{theorem}
\begin{remark}
	It is straightforward to repeat the above proof for the case of Dirichlet boundary conditions to obtain an analogous result for $\bigl\|\exp(-t(A^{\rm D}-\text{Id}))\bigr\|$. Here, the selfadjointness of $A^{\rm D}$ allows us to choose the half-angle $\theta$ arbitrarily close to $\pi/2$. 
\end{remark}

\section{Conclusion}\label{conclusion}

We have shown norm-resolvent convergence in the classical perforated domain problem with Dirichlet boundary conditions which has the interesting implication of spectral convergence (Cor. \ref{SpCon}). Some questions remain open and will be addressed in the future. 
 While the norm $\|J_\eps A_\eps^{-1}-A^{-1}J_\eps \|_{\mathcal L(L^2(\Om_\eps), L^2(\Om))}$ converges to 0, it is not clear from our method of proof how fast this happens. It would be desirable to obtain a precise convergence rate. In the case of Dirichlet boundary conditions a explicit convergence rate has been found by \cite{KP}.
 Another interesting question is whether in the case $\Om=\R^d$ there exist gaps in the spectrum of $A_\eps $ and how these depend on $\eps $. The existence of spectral gaps has been confirmed in two dimensions \cite{NRT}, but to the authors' knowledge the higher-dimensional case is still open.

\end{document}